\documentclass[12pt]{article}

\usepackage{amsmath,amsthm}
\usepackage{amssymb}
\usepackage{hyperref}

\def\sts#1#2{\left\{#1\atop#2\right\}}

\newtheorem{theorem}{Theorem}
\newtheorem{Prop}{Proposition}
\newtheorem{Cor}{Corollary}
\newtheorem{Lem}{Lemma}

\begin{document}

\title{Incomplete poly-Bernoulli numbers associated with incomplete Stirling numbers}
\author{
Takao Komatsu\footnote{The research of Takao Komatsu was supported in part by the grant of Wuhan University and by the Hubei Provincial Program}\\
\small School of Mathematics and Statistics\\[-0.8ex]
\small Wuhan University\\[-0.8ex]
\small Wuhan 430072 China\\[-0.8ex]
\small \texttt{komatsu@wuh.edu.cn}\\\\ 
K\'alm\'an Liptai\\ 
\small Rector\\[-0.8ex]
\small Eszterh\'azy K\'aroly College\\[-0.8ex] 
\small Eger 3300 Hungary\\[-0.8ex]
\small \texttt{liptaik@ektf.hu}\\\\
Istv\'an Mez\H{o}\footnote{The research of Istv\'an Mez\H{o} was supported by the Scientific Research Foundation of Nanjing University of Information Science \& Technology, and The Startup Foundation for Introducing Talent of NUIST. Project no.: S8113062001}\\
\small Department of Mathematics\\[-0.8ex]
\small Nanjing University of Information Science and Technology\\[-0.8ex]
\small Nanjing 210044 China\\[-0.8ex]
\small \texttt{istvanmezo81@gmail.com}
}

\date{
%\small Submitted: December 1, 2014;  Accepted: January 2, 2015.\\
%\small MR Subject Classifications: Primary 11B73; Secondary 05A18
}

\maketitle

\begin{abstract}
By using the associated and restricted Stirling numbers of the second kind, we give some generalizations of the poly-Bernoulli numbers. We also study their analytic and combinatorial properties. As an application, at the end of the paper we present a new infinite series representation of the Riemann zeta function via the Lambert $W$.
\end{abstract}

\section{Introduction}

Let $\mu\ge 1$ be an integer in the whole text. 
Our goal is to generalize the following relation for the poly-Bernoulli numbers $B_n^{(\mu)}$ (\cite[Theorem 1]{Kaneko}): 
$$
B_n^{(\mu)}=\sum_{k=0}^n(-1)^{n-k}\frac{k!}{(k+1)^\mu}\sts{n}{k}\quad(n\ge0,~\mu\ge 1)\,,
$$ 
where $\sts{n}{k}$ are the Stirling numbers of the second kind, determined by 
$$ 
\sts{n}{k}=\frac{1}{k!}\sum_{j=0}^k(-1)^i\binom{k}{j}(k-j)^n 
$$   
(see e.g., \cite{Jordan}).  
When $\mu=1$, $B_n^{(1)}$ are the classical Bernoulli numbers, defined by the generating function 
\begin{equation} 
\frac{x}{1-e^{-x}}=\sum_{n=0}^\infty B_n^{(1)}\frac{x^n}{n!}\,.
\label{gen:ber}
\end{equation}  
Notice that the classical Bernoulli numbers $B_n$ are also defined by the generating function 
$$
\frac{x}{e^x-1}=\sum_{n=0}^\infty B_n\frac{x^n}{n!}\,,
$$
satisfying $B_n^{(1)}=B_n$ ($n\ne 1$) with $B_1^{(1)}=1/2=-B_1$.  

The generating function of the poly-Bernoulli numbers $B_n^{(\mu)}$ is given by  
\begin{equation} 
\frac{{\rm Li}_\mu\bigl(1-e^{-x}\bigr)}{1-e^{-x}}=\sum_{n=0}^\infty B_n^{(\mu)}\frac{x^n}{n!}\,, 
\label{gen:pber}
\end{equation} 
where 
$$
{\rm Li}_\mu(z)=\sum_{m=1}^\infty\frac{z^m}{m^\mu}
$$
is the $\mu$-th polylogarithm function (\cite[(1)]{Kaneko}). The generating function of the poly-Bernoulli numbers can also be written in terms of iterated integrals (\cite[(2)]{Kaneko}): 
\begin{equation} 
e^x\cdot\underbrace{\frac{1}{e^x-1}\int_0^x\frac{1}{e^x-1}\int_0^x\dots\frac{1}{e^x-1}\int_0^x}_{\mu-1}\frac{x}{e^x-1}\underbrace{dx dx\dots dx}_{\mu-1}
=\sum_{n=0}^\infty B_n^{(\mu)}\frac{x^n}{n!}\,. 
\label{gen:iter}
\end{equation}  

Several generalizations of the poly-Bernoulli numbers have been considered (\cite{BH1,BH2,CC,Jolany,Sasaki}). However, most kinds of generalizations are based upon the generating functions of (\ref{gen:ber}) and/or (\ref{gen:pber}). On the contrary, our generalizations are based upon the explicit formula in terms of the Stirling numbers. In \cite{KMS}, a similar approach is used to generalize the Cauchy numbers $c_n$, defined by $x/\log(1+x)=\sum_{n=0}^\infty c_n x^n/n!$. In this paper, by using the associated and restricted Stirling numbers of the second kind, we give substantial generalizations of the poly-Bernoulli numbers. One of the main results is to generalize the formula in (\ref{gen:pber}) as 
$$ 
\sum_{n=0}^\infty B_{n,\le m}^{(\mu)}\frac{t^n}{n!}=\frac{{\rm Li}_\mu\bigl(1-E_m(-t)\bigr)}{1-E_m(-t)}
$$ 
and
$$  
\sum_{n=0}^\infty B_{n,\ge m}^{(\mu)}\frac{t^n}{n!}=\frac{{\rm Li}_\mu\bigl(E_{m-1}(-t)-e^{-t}\bigr)}{E_{m-1}(-t)-e^{-t}}\,, 
$$   
where $E_m(t)=\sum_{k=0}^m\frac{t^k}{k!}$.  
See Theorem \ref{th10} below.

\section{Incomplete Stirling numbers of the second kind}

In place of the classical Stirling numbers of the second kind $\sts{n}{k}$ we substitute the restricted Stirling numbers and the associated Stirling numbers
$$
\sts{n}{k}_{\le m}\quad\mbox{and}\quad\sts{n}{k}_{\ge m}\,,
$$
respectively.  
Some combinatorial and modular properties of these numbers can be found in \cite{Mezo}, and other properties can be found in the cited papers of \cite{Mezo}.
The generating functions of these numbers are given by 
\begin{equation} 
\sum_{n=k}^{mk}\sts{n}{k}_{\le m}\frac{x^n}{n!}=\frac{1}{k!}(E_m(x)-1)^k
\label{gen:reststs}
\end{equation} 
and
\begin{equation} 
\sum_{n=mk}^\infty\sts{n}{k}_{\ge m}\frac{x^n}{n!}=\frac{1}{k!}\bigl(e^x-E_{m-1}(x)\bigr)^k
\label{gen:assosts} 
\end{equation} 
respectively, 
where 
$$ 
E_m(t)=\sum_{k=0}^m\frac{t^k}{k!}
$$ 
is the $m$th partial sum of the exponential function sum. 
These give the number of the $k$-partitions of an $n$-element set, such that each block contains at most or at least $m$ elements, respectively. 
Since the generating function of $\sts{n}{k}$ is given by 
$$
\sum_{n=k}^\infty\sts{n}{k}\frac{x^n}{n!}=\frac{(e^x-1)^k}{k!}
$$ 
(see e.g., \cite{Jordan}), 
by $E_\infty(x)=e^x$ and $E_0(x)=1$, we have 
$$
\sts{n}{k}_{\le\infty}=\sts{n}{k}_{\ge 1}=\sts{n}{k}\,. 
$$ 
These give the number of the $k$-partitions of an $n$-element set, such that each block contains at most or at least $m$ elements, respectively.
Notice that these numbers where $m=2$ have been considered by several authors (e.g., \cite{Comtet, Howard1, Rio, Zhao}).  
\medskip

It is well-known that the Stirling numbers of the second kind satisfy the recurrence relation: 
\begin{equation} 
\sts{n+1}{k}=k\sts{n}{k}+\sts{n}{k-1}
\label{rec:sts} 
\end{equation} 
for $k>0$, with the initial conditions
$$
\sts{0}{0}=1\quad\mbox{and}\quad \sts{n}{0}=\sts{0}{n}=0
$$
for $n>0$. 
The restricted and associated Stirling numbers of the second kind satisfy the similar relations.  It is easy to see the initial conditions
\begin{align*}
&\sts{0}{0}_{\le m}=1\quad\mbox{and}\quad \sts{n}{0}_{\le m}=\sts{0}{n}_{\le m}=0\,,\\
&\sts{0}{0}_{\ge m}=1\quad\mbox{and}\quad \sts{n}{0}_{\ge m}=\sts{0}{n}_{\ge m}=0
\end{align*}
for $n>0$.

\begin{Prop}
For $k>0$ we have 
\begin{align} 
\sts{n+1}{k}_{\le m}&=\sum_{i=0}^{m-1}\binom{n}{i}\sts{n-i}{k-1}_{\le m}
\label{rec:rests}\\
&=k\sts{n}{k}_{\le m}+\sts{n}{k-1}_{\le m}-\binom{n}{m}\sts{n-m}{k-1}_{\le m}\,,
\label{rec:rests2}\\
\sts{n+1}{k}_{\ge m}&=\sum_{i=m-1}^n\binom{n}{i}\sts{n-i}{k-1}_{\ge m}
\label{rec:assos}\\
&=k\sts{n}{k}_{\ge m}+\binom{n}{m-1}\sts{n-m+1}{k-1}_{\ge m}\,.
\label{rec:assos2}
\end{align} 
\label{prp1}
\end{Prop} 

\noindent 
{\it Remark.}  
The fourth relation (\ref{rec:assos2}) appeared in a different form in \cite{Howard1}. 
Since 
$$
\sum_{i=1}^{n-k+1}\binom{n}{i}\sts{n-i}{k-1}=k\sts{n}{k}\,,  
$$ 
the relation (\ref{rec:rests}) and the relation (\ref{rec:assos}) are both reduced to the relation (\ref{rec:sts}), if $m\ge n-k+2$ and if $m=1$, respectively. It is trivial to see that the relations (\ref{rec:rests2}) and (\ref{rec:assos2}) are also reduced to the original relation (\ref{rec:sts}), if $m>n$ and if $m=1$, respectively.

\noindent 
{\it Proof of Proposition \ref{prp1}.}  
The combinatorial proofs of the previous theorem are given as follows.
We shall give combinatorial proofs.  First, 
identity \eqref{rec:rests}. To construct a partition with $k$ blocks on $n+1$ element we can do the following. The last element in its block can have $i$ elements by side, where $i=0,1,\dots,m-1$. We have to choose these $i$ elements from $n$. This can be done in $\binom{n}{i}$ ways. The rest of the elements go into $k-1$ blocks in $\sts{n-i}{k-1}_{\le m}$ ways. Summing over the possible values of $i$ we are done.

The proof of \eqref{rec:rests2}. The above construction can be described in another way: the last element we put into a singleton and the other $n$ elements must form a partition with $k-1$ blocks: $\sts{n}{k-1}_{\le m}$ possibilities. Or we put this element into one existing block after constructing a partition of $n$ elements into $k$ blocks. This offers us $k\sts{n}{k}_{\le m}$ possibilities, but we must subtract the possibilities when we exceed the block size limit $m$. This happens if we put the last element into a block of $m$ elements. There are $\binom{n}{m}\sts{n-m}{k-1}_{\le m}$ such partitions in total. The proof is done.

The proof of \eqref{rec:assos} and \eqref{rec:assos2} is similar.
\qed

\medskip 

Note that the classical Stirling numbers of the second kind $\sts{n}{k}$ satisfy the identities:  
\begin{align*}
%&\sts{0}{0}=1,\quad \sts{n}{0}=0\quad(n>0),\quad \sts{0}{k}=0\quad(k>0),\quad \sts{n}{k}=0\quad(k>n),\\
&\sts{n}{1}=\sts{n}{n}=1,\quad \sts{n}{n-1}=\binom{n}{2},\\
&\sts{n}{n-2}=\frac{3n-5}{4}\binom{n}{3},\quad \sts{n}{n-3}=\binom{n}{4}\binom{n-2}{2},\\
&\sts{n}{2}=2^{n-1}-1,\quad \sts{n}{3}=\frac{3^{n-1}}{2}-2^{n-1}+\frac{1}{2},\quad \\
&\sts{n}{4}=\frac{4^{n-1}}{6}-\frac{3^{n-1}}{2}+2^{n-2}-\frac{1}{6}\,. 
\end{align*}

By the definition (\ref{gen:reststs}) or Proposition \ref{prp1} (\ref{rec:rests}), we list several basic properties about the restricted Stirling numbers of the second kind.  Some basic properties about the associated Stirling numbers of the second kind can be found in \cite{Comtet,Howard1,Mezo,Zhao}.  

\begin{Lem}  
For $0\le n\le k-1$ or $n\ge m k+1$, we have 
\begin{equation}
\sts{n}{k}_{\le m}=0\,. 
\label{le00}
\end{equation} 
For $k\le n\le m k$, we have 
\begin{align}
\sts{n}{k}_{\le m}&=\sts{n}{k}\quad(k\le n\le m)\,,
\label{le01}\\
\sts{n}{n}_{\le m}&=1\quad(n\ge 0,~m\ge 1)\,,
\label{le02}\\
\sts{n}{n-1}_{\le m}&=\binom{n}{2}\quad(n\ge 2,~m\ge 2)\,,
\label{le03}\\
\sts{n}{n-2}_{\le m}&=
\begin{cases} 
\frac{3n-5}{4}\binom{n}{3}&(n\ge 4,~m\ge 3);\\
3\binom{n}{4}&(n\ge 4,~m= 2)\,,
\end{cases} 
\label{le04}\\ 
\sts{n}{n-3}_{\le m}&=
\begin{cases} 
\binom{n}{4}\binom{n-2}{2}&(n\ge 4,~m\ge 4);\\
15\binom{n}{6}+10\binom{n}{5}&(n\ge 4,~m=3);\\
15\binom{n}{6}&(n\ge 4,~m=2)\,,
\end{cases} 
\label{le05}\\ 
\sts{n}{1}_{\le m}&=1\quad(1\le n\le m)\,,
\label{le06}\\
\sts{n}{2}_{\le m}&=2^{n-1}-1\quad(2\le n\le m+1)\,.
\label{le07}
\end{align} 
\label{lemreststs}  
\end{Lem}  
\begin{proof}  
Some of the above special values are trivial. Some of them can be proven by analyzing the possible block structures.

We take \eqref{le04} as a concrete example.

Let $m=2$, and the number of blocks be $k=n-2$. Then for the block structure we have the only one possibility
\[\underbrace{\,.\,|\,.\,|\cdots|\,.\,|}_{n-4}\,.\,.\,|\,.\,.\]
That is, there are $n-4$ singletons and two blocks of length 2. There are $\frac12\binom{4}{2}\binom{n}{4}=3\binom{n}{4}$ such partitions: we have to choose those four elements going to the non singleton blocks in $\binom{n}{4}$ ways. Then we put two of four into the first block and the other two goes to the other: $\binom{4}{2}=6$ cases. Finally, we have to divide by two because the order of the blocks does not matter. The last case of \eqref{le04} follows.

If $m=3$ then we have one more possible distribution of blocks sizes apart from the above:
\[\underbrace{\,.\,|\,.\,|\cdots|\,.\,|}_{n-5}\,.\,.\,.\,\]
Into the last block we have $\binom{n}{3}$ possible option to put 3 elements. So if $m=3$ and $k=n-2$ then we have $\binom{n}{3}+3\binom{n}{4}=\frac{3n-5}{4}\binom{n}{3}$ cases in total.

The rest of the cases can be treated similarly.
\end{proof}

\section{Incomplete poly-Bernoulli numbers}

\subsection{Generating function and its integral representation}

By using two types of incomplete Stirling numbers, define {\it restricted poly-Bernoulli numbers} $B_{n,\le m}^{(\mu)}$ and {\it associated poly-Bernoulli numbers} $B_{n,\ge m}^{(\mu)}$ by 
\begin{equation} 
B_{n,\le m}^{(\mu)}=\sum_{k=0}^n(-1)^{n-k}\frac{k!}{(k+1)^\mu}\sts{n}{k}_{\le m}\quad(n\ge0),
\label{def:assobern}
\end{equation} 
and
\begin{equation} 
B_{n,\ge m}^{(\mu)}=\sum_{k=0}^n(-1)^{n-k}\frac{k!}{(k+1)^\mu}\sts{n}{k}_{\ge m}\quad(n\ge0)\,,
\label{def:restbern}
\end{equation} 
respectively.  
These numbers can be considered as generalizations of the usual poly-Bernoulli numbers $B_n^{(\mu)}$, since
$$
B_{n,\le\infty}^{(\mu)}=B_{n,\ge 1}^{(\mu)}=B_n^{(\mu)}\,.
$$ 
We call these numbers as {\it incomplete poly-Bernoulli numbers}. 
\medskip

One can deduce that these numbers have the generating functions. 

\begin{theorem} 
We have 
\begin{equation} 
\sum_{n=0}^\infty B_{n,\le m}^{(\mu)}\frac{t^n}{n!}=\frac{{\rm Li}_\mu\bigl(1-E_m(-t)\bigr)}{1-E_m(-t)}
\label{gen:restbern}
\end{equation}
and
\begin{equation} 
\sum_{n=0}^\infty B_{n,\ge m}^{(\mu)}\frac{t^n}{n!}=\frac{{\rm Li}_\mu\bigl(E_{m-1}(-t)-e^{-t}\bigr)}{E_{m-1}(-t)-e^{-t}}\,.
\label{gen:assobern}
\end{equation}  
\label{th10}  
\end{theorem}  

\noindent 
{\it Remark.} 
In the first formula $m\to\infty$ gives back the poly-Bernoulli numbers (\ref{gen:pber}) since $E_\infty(-t)=e^{-t}$ and ${\rm Li}_1(z)=-\log(1-z)$, while in the second we must take $m=1$ since $E_0(-t)=1$.  

\noindent
{\it Proof of Theorem \ref{th10}.} 
By the definition of (\ref{def:restbern}) and using (\ref{gen:reststs}), we get 
\begin{align*} 
\sum_{n=0}^\infty B_{n,\le m}^{(\mu)}\frac{t^n}{n!}&=\sum_{n=0}^\infty\sum_{k=0}^n\sts{n}{k}_{\le m}\frac{(-1)^{n-k}k!}{(k+1)^\mu}\frac{t^n}{n!}\\
&=\sum_{k=0}^\infty\frac{(-1)^k k!}{(k+1)^\mu}\sum_{n=k}^\infty\sts{n}{k}_{\le m}\frac{(-t)^n}{n!}\\
&=\sum_{k=0}^\infty\frac{(-1)^k k!}{(k+1)^\mu}\frac{1}{k!}\left((-t)+\frac{(-t)^2}{2!}+\cdots+\frac{(-t)^m}{m!}\right)^k\\
&=\sum_{k=0}^\infty\frac{\bigl(1-E_m(-t)\bigr)^k}{(k+1)^\mu}\\
&=\frac{{\rm Li}_\mu\bigl(1-E_m(-t)\bigr)}{1-E_m(-t)}\,.
\end{align*} 
Similarly, 
by the definition of (\ref{def:assobern}) and using (\ref{gen:assosts}), we get 
\begin{align*} 
\sum_{n=0}^\infty B_{n,\ge m}^{(\mu)}\frac{t^n}{n!}&=\sum_{n=0}^\infty\sum_{k=0}^n\sts{n}{k}_{\le m}\frac{(-1)^{n-k}k!}{(k+1)^\mu}\frac{t^n}{n!}\\
&=\sum_{k=0}^\infty\frac{(-1)^k k!}{(k+1)^\mu}\sum_{n=k}^\infty\sts{n}{k}_{\ge m}\frac{(-t)^n}{n!}\\
&=\sum_{k=0}^\infty\frac{(-1)^k k!}{(k+1)^\mu}\frac{1}{k!}\left(\frac{(-t)^m}{m!}+\frac{(-t)^{m+1}}{(m+1)!}+\cdots\right)^k\\
&=\sum_{k=0}^\infty\frac{\bigl(E_{m-1}(-t)-e^{-t}\bigr)^k}{(k+1)^\mu}\\
&=\frac{{\rm Li}_\mu\bigl(E_{m-1}(-t)-e^{-t}\bigr)}{E_{m-1}(-t)-e^{-t}}\,.
\end{align*} 
\qed

For $\mu\ge 1$, the generating functions can be written in the form of iterated integrals.  We set $E_{-1}(-t)=0$ for convenience. 

\begin{theorem}
\begin{align}
&\frac{1}{1-E_m(-t)}\cdot\underbrace{\int_0^t\frac{E_{m-1}(-t)}{1-E_m(-t)}\int_0^t\dots\frac{E_{m-1}(-t)}{1-E_m(-t)}\int_0^t\frac{E_{m-1}(-t)}{1-E_m(-t)}}_{\mu-1}\notag\\
&\qquad\times\bigl(-\log\bigl(E_m(-t)\bigr)\bigr)\underbrace{dt\dots dt}_{\mu-1}\notag\\
&\qquad\qquad =\sum_{n=0}^\infty B_{n,\le m}^{(\mu)}\frac{x^n}{n!}
\label{gen:iterrest}\\
&\frac{1}{E_{m-1}(-t)-e^{-t}}\cdot\underbrace{\int_0^t\frac{e^{-t}-E_{m-2}(-t)}{E_{m-1}(-t)-e^{-t}}\int_0^t\dots\frac{e^{-t}-E_{m-2}(-t)}{E_{m-1}(-t)-e^{-t}}\int_0^t\frac{e^{-t}-E_{m-2}(-t)}{E_{m-1}(-t)-e^{-t}}}_{\mu-1}\notag\\
&\qquad\times\bigl(-\log\bigl(1+e^{-t}-E_{m-1}(-t)\bigr)\bigr)\underbrace{dt\dots dt}_{\mu-1}\notag\\
&\qquad\qquad =\sum_{n=0}^\infty B_{n,\ge m}^{(\mu)}\frac{x^n}{n!}\,. 
\label{gen:iterasso}
\end{align}
\label{th20} 
\end{theorem}  

\noindent 
{\it Remark.}  
If $m\to\infty$ in (\ref{gen:iterrest}), by $E_\infty(-t)=e^{-t}$, and if $m=1$ in (\ref{gen:iterasso}), by $E_0(-t)=1$ and $E_{-1}(-t)=0$, both of them are reduced to (\ref{gen:iter}). 

\noindent 
{\it Proof of Theorem \ref{th20}.} 
Since for $\mu\ge 1$ 
$$
\frac{d}{dt}{\rm Li}_\mu\bigl(1-E_m(-t)\bigr)=\frac{E_{m-1}(-t)}{1-E_m(-t)}{\rm Li}_{\mu-1}\bigl(1-E_m(-t)\bigr)\,,
$$ 
we have
\begin{align*} 
&{\rm Li}_\mu\bigl(1-E_m(-t)\bigr)\\
&=\int_0^t\frac{E_{m-1}(-t)}{1-E_m(-t)}{\rm Li}_{\mu-1}\bigl(1-E_m(-t)\bigr)dt\\
&=\int_0^t\frac{E_{m-1}(-t)}{1-E_m(-t)}\int_0^t\frac{E_{m-1}(-t)}{1-E_m(-t)}{\rm Li}_{\mu-2}\bigl(1-E_m(-t)\bigr)dt dt\\
&=\int_0^t\frac{E_{m-1}(-t)}{1-E_m(-t)}\int_0^t\frac{E_{m-1}(-t)}{1-E_m(-t)}\cdots\int_0^t\frac{E_{m-1}(-t)}{1-E_m(-t)}{\rm Li}_1\bigl(1-E_m(-t)\bigr)\underbrace{dt\cdots dt}_{\mu-1}\\
&=\int_0^t\frac{E_{m-1}(-t)}{1-E_m(-t)}\int_0^t\frac{E_{m-1}(-t)}{1-E_m(-t)}\cdots\int_0^t\frac{E_{m-1}(-t)}{1-E_m(-t)}\bigl(-\log\bigl(E_m(-t)\bigr)\bigr)\underbrace{dt\cdots dt}_{\mu-1}\,.
\end{align*} 
Therefore, we obtain (\ref{gen:iterrest}). 
Similarly, by 
$$
\frac{d}{dt}{\rm Li}_\mu\bigl(E_{m-1}(-t)-e^{-t}\bigr)=\frac{e^{-t}-E_{m-2}(-t)}{E_{m-1}(-t)-e^{-t}}{\rm Li}_{\mu-1}\bigl(E_{m-1}(-t)-e^{-t}\bigr)\,,
$$ 
we obtain (\ref{gen:iterasso}).
\qed

If $\mu=1$ in Theorem \ref{th10} or in Theorem \ref{th20}, the generating functions of the {\it restricted Bernoulli numbers} $B_{n,\le m}^{(1)}$ and {\it associated Bernoulli numbers} $B_{n,\ge m}^{(1)}$ are given.  Both functions below are reduced to the generating function (\ref{gen:ber}) of the Bernoulli numbers $B_n^{(1)}$ if $m\to\infty$ and $m=1$, respectively. 

\begin{Cor}  
We have 
$$ 
\sum_{n=0}^\infty B_{n,\le m}^{(1)}\frac{t^n}{n!}=\frac{\log E_m(-t)}{E_m(-t)-1}
$$
and
$$ 
\sum_{n=0}^\infty B_{n,\ge m}^{(1)}\frac{t^n}{n!}=\frac{\log\bigl(1+e^{-t}-E_{m-1}(-t)\bigr)}{e^{-t}-E_{m-1}(-t)}\,.
$$   
\label{cor11} 
\end{Cor}

\subsection{Basic divisibility for non-positive $\mu$}

In this short subsection we deduce a basic divisibility property for both the restricted and associated poly-Bernoulli numbers.

It is known \cite{GKP} that
\[\sts{p}{k}\equiv0\pmod{p}\quad(1<k<p)\]
for any prime $p$. The proof of this basic divisibility is the same for the restricted and associated Stirling numbers, so we can state that
\[\sts{p}{k}_{\le m}\equiv0\pmod{p}\quad(k=0,1,\dots),\]
and
\[\sts{p}{k}_{\ge m}\equiv0\pmod{p}\quad(k\ge2).\]
(Note that $\sts{p}{1}_{\ge m}=1$.) These immediately lead to the next statement.

\begin{theorem}For any $\mu\le0$ we have that
\begin{align*}
B_{p,\le m}^{(\mu)}\equiv&0\pmod{p},\\
B_{p,\ge m}^{(\mu)}\equiv&2^{|\mu|}\pmod{p}
\end{align*}
hold for any prime $p$.
\end{theorem}

\section{A new series representation for the Riemann zeta function}

To present our result, we need to recall the definition of the Lambert $W$ function. $W(a)$ is the solution of the equation
\[xe^x=a,\]
that is, $W(a)e^{W(a)}=a$. Since this equation, in general, has infinitely many solutions, the $W$ function has infinitely many complex branches denoted by $W_k(a)$ where $k\in\mathbb Z$. What we prove is the following.

\begin{theorem}For any $\mu\in\mathbb C$ with $\Re(\mu)>1$ we have that
\[\zeta(\mu)=\sum_{n=0}^\infty B_{n,\ge2}^{(\mu)}\frac{(W_k(-1))^n}{n!}\]
for $k=0,-1$, where $\zeta$ is the Riemann zeta function.
\end{theorem}

\begin{proof} Let us recall the generating function of $B_{n,\ge m}^{(\mu)}$ in the particular case when $m=2$:
\begin{equation}
\sum_{n=0}^\infty B_{n,\ge2}^{(\mu)}\frac{(-t)^n}{n!}=\frac{{\rm Li}_{\mu}(1+t-e^t)}{1+t-e^t}.\label{zetapr}
\end{equation}
By a simple transformation it can be seen that the equation $1+t-e^t=1$ is solvable in terms of the Lambert $W$ function, and that the solution is $-W_k(-1)$ for any branch $k\in\mathbb Z$. However, \eqref{zetapr} is valid only for $t$ such that $|1+t-e^t|\le 1$, at least when $\Re(\mu)>1$. (This comes from the proof of Theorem \ref{th10}.) Since the absolute value of $-W_k(-1)$ grows with $k$, the only two branches which belong the convergence domain of \eqref{zetapr} is $k=-1,0$. Hence, substituting one of these in place of $t$ we have that
\[\sum_{n=0}^\infty B_{n,\ge2}^{(\mu)}\frac{(W_k(-1))^n}{n!}=\frac{{\rm Li}_{\mu}(1-W_k(-1)-e^{-W_k(-1)})}{1-W_k(-1)-e^{-W_k(-1)}}=\frac{{\rm Li}_{\mu}(1)}{1}=\zeta(\mu).\]
\end{proof}
Note that
\[W_{0}(-1)=\overline{W_{-1}(-1)}\approx-0.318132+1.33724 i,\]
so all the terms in the incomplete Bernoulli sum are complex, but the sum itself always converges to the real number $\zeta(\mu)$.

\section{Acknowledgement}  
This work has been partly done when the first author stayed in Eszterh\'azy K\'aroly College by Balassi Institute Program in 2014.

%\noindent 2010 {\it Mathematics Subject Classification}:
%Primary 11B73; Secondary 05A18

%\noindent \emph{Keywords:} 
%Stirling numbers, poly-Bernoulli numbers, Bernoulli numbers

\end{document}